\documentclass[]{article}
\usepackage{indentfirst}
\usepackage{cite}
\usepackage{amsfonts}
\usepackage{amsmath}

\title{\textbf{Area bounds for minimal surfaces in geodesic ball of hyperbolic space} }
\date{}
\author{Jingze Zhu\footnote{Tsinghua University, department of Mathematical Sciences. This article was done during my one-month visit to Columbia University. I thank to Professor Mu-Tao Wang for providing the suggestion to study this kind of generalization.}}

\begin{document}

\maketitle
\begin{abstract}
  In hyperbolic space $H^n$ we set a geodesic ball of radius $\rho$. Consider a $k$ dimensional minimal submanifold passing through the origin of the geodesic ball with boundary lies on the boundary of that geodesic ball. We prove that its area is no less than the totally geodesic $k$ dimensional submanifold passing through the origin in that geodesic ball. This is a partial generalization of the corresponding problem in $R^n$.
\end{abstract}
\section{Introduction}
This article mainly proves a sharp lower bound of area for $k$ dimensional minimal surfaces passing through the origin of the geodesic ball $B_p(\rho_0) $ in standard hyperbolic space $H^n$ with boundary on $ \partial B_p(\rho_0) $

Fix a point $p \in H^n$ and a positive number $\rho_0$. Denote the geodesic ball in $ H^n $ with center $p$ and radius $\rho_0$ by $B_p(\rho_0)$.
Consider the $k$ dimensional plane $P$ passing through the origin in $T_p H^n$. Let $N_0=\textnormal{exp}_p(P) \cap B_p(\rho_0)$. Then $N_0$ is a totally geodesic surface and thus a minimal surface in $H^n$ since its second fundamental form vanishes.

Now we can state the theorem:

\newtheorem{thm}{Theorem}
\newtheorem{lem}{Lemma}
\begin{thm}
  Suppose $M$ is a $k$ dimensional minimal surface of $H^n$ passing through $p$ with the boundary $\partial M$ lies on
  $\partial B_p(\rho_0)$. Then the area of $M$ is no less than the area of $N_0$, i.e.  $|M|\geq|N_0|$
\end{thm}

It is one generalization of a result on the area bound of a minimal surface passing through the a fixed point $P$ of a ball in $R^n$ with boundary on the ball. If $P=O$ be the origin of the ball, then by monotonicity formula the area is no less than a plane. It has also been proved recently in \cite{simonbrendlepeikenhung.{2016}} that if $P\neq O$, the area is no less than the area of the plane passing through $P$, perpendicular to $OP$. The technique it used was to find a proper vector field and use the first variation formula. Inspired from this, we can proceed to Theorem 1. In fact, the vector field we finally derived looks like the gradient of the Green function of $B_p(\rho_0)$ (with power changed).

\section{Proof of the theorem}
Let $|N_0|=\omega$ . We shall find a vector field $W$ defined on $H^n \backslash \{p\}$ that satisfies following three conditions
\begin{enumerate}
  \item $W$ vanishes on $\partial B_p(\rho_0)$
  \item $|W(x)| \rightarrow -\frac{\omega r^{-{k+1}}}{\omega_{k-1}}$ as $x \rightarrow p$

  where $r=r(x,p)$ is the distance between $p$ and $x$ and $\omega_{k-1}$ is area of $k-1$ dimensional unit sphere.
  \item $ \sum_{i=1}^{k} \langle D_{\tau _i}W, \tau_i \rangle \leq 1 $

  where $D$ is the covariant derivative for $H^n$ and $\{\tau_i , i=1,2,...,k\} $ are $k$ orthonormal vectors in $T_qH^n$
  for any point $q \in H^n \backslash \{p\} $
\end{enumerate}
 We now turn our attention to the rotational symmetric vector field. Choose Gauss Normal Coordinate of point $p$. More explicitly,
 Let $\{\Theta \in S^{n-1} , r \}$ be canonical polar coordinate for $T_pH^n$. Then using exponential map to get the coordinate for $H^n\backslash \{p\}$.

 By Gauss's theorem,  $\frac{\partial}{\partial r} = \nabla r $ is unit vector perpendicular to $\partial B_p(r)$ . For simplicity, let
 $\partial_r = \frac{\partial}{\partial r}$.

 Let $f=f(r)$ be a smooth function defined on $H^n$, depending only on the distance to $p$. Let $W=f(r)\partial_r$. We have the following lemma.

\begin{lem}
   Let $q\in H^n \neq p$. For orthonormal vectors $\{\tau_i, i=1,2,...,k\} \in T_qH^n$, r(q)=r(q,p), we have
 $\sum_{i=1}^{k} \langle D_{\tau_i}W, \tau_i \rangle = kf \coth r + (f'-f \coth r) \cdot |(\partial_r)^T|^2$,
  where $(\partial_r)^T$ is orthogonal projection of $\partial r$ onto the vector space spanned by $\{\tau_i\}$, i.e.
   $(\partial_r)^T = \langle \partial_r , \tau_i \rangle \tau_i  $
\end{lem}
   Proof: Using the metric form $ds^2=dr^2+(\sinh r)^2 \ d\Omega^2$ for $H^n$ where $d \Omega ^2$ is the standard metric for unit sphere $S^{n-1}$.

   Fix an $r=r_0$.  For any point $q \neq p$ in $H^n$, let $\{e_i,i=1,2,...,n-1\}$ be orthonormal basis for $T_q\partial B_p(r)$ and let $e_n=\partial_r$. Since the induced metric on $\partial B_p(r_0)$, i.e. the level set of $r$, is spherical symmetric metric depending only on $r$, by symmetry we have $D_{e_i}\partial _r = \psi e_i, i=1,...,n-1$ where $\psi =\psi(r)$ is a function depending only on $r$, i.e $D \partial_r$ as a linear transformation is a multiple of identity.
   Further, $ D_{\partial_r}\partial_r =0$. Hence the Laplacian $\Delta r = (n-1)\psi(r)$.

   However we have already known that in a space of constant negative sectional curvature $-\lambda ^2$  with $r$ be distance function,  $\Delta r = (n-1)\lambda \coth  \lambda r $. Let $\lambda=1$, we will have $\psi(r)=\coth r$.

    Notice that $e_i(f)=0$ for $i=1,...,n-1$ since $f$ is a constant on $\partial B_p(r)$.Then

   \[D_{e_i}W= f(r)D_{e_i}e_n = f (\coth r) e_i \ \ \ \  i=1,2,...,n-1\]

   \[D_{e_n}W=D_{e_n}f(r)e_n + fD_{e_n}e_n=f'(r)e_n \]

   Let $\{\tau_j, j=1,2,...,k\}$be $k$ orthonormal vectors in $T_qH^n$.

   Assume $\tau_j=\sum_{i=1}^{n}a_{ji}e_i$. Thus $A=(a_{ji})_{k \times n}$ is a matrix with $k$ orthonormal row vectors.

   After doing that, we will have

   \[\sum_{j=1}^{k}\langle D_{\tau_j}W,\tau_j \rangle = kf \coth r + (f'-f \coth r)(\sum_{j=1}^{k}a_{jn}^2) \]

   On the other hand

    \[(\sum_{j=1}^{k}a_{jn}^2) = \sum_{j=1}^{k}\langle e_n, \tau_j \rangle ^2 = |\sum_{j=1}^{k}\langle e_n, \tau_j \rangle \tau_j|^2=|(\partial_r)^T|^2\]

   So we have proven Lemma1.

   \begin{lem}
   Let \[f(r)=(\sinh r)^{-(k-1)} \int_{\rho_0}^{r}(\sinh \rho)^{k-1} \mathrm d \rho\] then vector field $W=f(r)\partial_r$\ defined on $B_p(\rho_0) \backslash \{p\}$ satisfies all the three conditions enumerated above.
\end{lem}
Proof: Since $f(\rho_0)=0$, then $W$ vanishes on $\partial B_p(\rho_0)$.

 We can write \[f(r)=-C(\sinh r)^{-(k-1)}+(\sinh r)^{-(k-1)} \int_{0}^{r}(\sinh \rho)^{k-1} \mathrm d \rho \]

where \[C= \int_{0}^{\rho_0}(\sinh \rho)^{k-1} \mathrm d \rho\] is a constant depending only on $\rho_0$.

Then by direct computation we can check $f$ defined above solves the equation

 \[f'+(k-1)f \coth r =1\] with boundary condition $f(\rho_0)=0$.

 With the explicit form of the solution given above we can see $f\leq 0 $ on $(0,\rho_0]$. Thus \ $f'=1-(k-1)f \coth r >0$. \ Then $f'-f\coth r >0 $ on $(0,\rho_0]$.

 Since $|\partial_r|=1$, thus $|(\partial_r)^T| \leq 1$ and equality holds iff $\partial_r $ lies on the subspace spanned by ${\tau_1,...,\tau_k}$

 Together with lemma1 ,we have
\begin{eqnarray*}
\sum_{i=1}^{k} \langle D_{\tau_i}W, \tau_i \rangle &=& kf \coth r + (f'-f \coth r) \cdot |(\partial_r)^T|^2 \\
&\leq& kf \coth r + (f'-f \coth r) \cdot 1= f'+(k-1)f \coth r =1
\end{eqnarray*}

 Moreover, equality holds iff $\partial_r$ lies on the subspace spanned by $\tau_1,...,\tau_k$
 and condition3 is satisfied.

 Using the metric of the form $ds^2=dr^2+(\sinh r)^2 \ d\Omega^2 $ , we can see that area of a totally geodesic surface in $H^n$ is
 \[\omega=|N_0|=\int_{0}^{\rho_0}\omega_{k-1}(\sinh \rho)^{k-1} \mathrm d \rho\]

 But $\sinh r \rightarrow r $ as $r \rightarrow 0$, and \[(\sinh r)^{-(k-1)} \int_{0}^{r}(\sinh \rho)^{k-1} \mathrm d \rho = O(r) \rightarrow 0 \ \ \ \mathrm{as} \ \ \ r \rightarrow 0\]

 So we have \[f \rightarrow \frac{\omega r^{-{k+1}}}{\omega_{k-1}}\ \  \mathrm{as} \ \ r \rightarrow 0\]
 Thus \[|W(x)| \rightarrow \frac{\omega r^{-{k+1}}}{\omega_{k-1}} \ \ as \ \  x \rightarrow p\]
 Then condition2 is satisfied. \\ \\

 Now we prove the theorem with the help of the vector field $W$.

 Let $M$ be any minimal surface in $H^n$ passing through point $p$ and  with boundary on $\partial B_p(\rho_0)$. At each point $q \in M\backslash \{p\}$, let $\tau_1, ...,\tau_k$ \ be orthonormal basis for $T_qM$ and let $W^T$ be orthogonal projection of $W$ onto $T_qM$.

 Let $H$ be mean curvature vector about the immersion $M \subset H^n$ and $div_M$ is the divergence operator on submanifold $M$, it is known that
\[div_M W^T = \sum_{i=1}^{k} \langle D_{\tau_i}W, \tau_i \rangle + \langle H, W\rangle \]

Since $M$ is minimal and thus $H $ \ is zero vector, then

\[div_M W^T = \sum_{i=1}^{k} \langle D_{\tau_i}W, \tau_i \rangle \leq 1\] on $M\backslash  \{p\}$\ \ \ by property 3 of $W$.

Using divergence formula,for each sufficient small $r>0$ \ we have

\[|M| \geq \int_{M\backslash B_p(r)} div_M W^T = \int_{M \cap \partial B_p(\rho_0)}\langle W^T, \nu \rangle + \int_{M \cap \partial B_p(r)}\langle W^T, \nu \rangle \]
where $\nu$ is outside unit vector perpendicular to boundary of $M \backslash B_p(r)$(in the sense of submanifold).
Since $W$ vanishes on $\partial B_p(\rho_0)$, the first term on the  $\mathcal{RHS}$ of above formula vanishes.
Moreover,on $M \cap \partial B_p(r), \nu = -\partial_r + o(r)$.

 Together with asymptotic behaviour of $W$, we have

  \[|M|\geq \lim_{r \rightarrow 0}\int_{M \cap \partial B_p(r)}\langle W^T, \nu \rangle \] \[= \lim_{r \rightarrow 0}\int_{M \cap \partial B_p(r)} -\frac{\omega}{\omega_{k-1}}r^{-(k-1)}(-1+o(r)) \\ = \omega = |N_0|\] which proves the theorem. Obviously, choose $M=N_0$ will attain the minimum.

\section{A generalization}
 In fact we can derive the formula for the space of any space forms. What really matters is the area of a unit sphere and the coefficient $\phi(r)$ in warped product metric $ds^2=dr^2+\phi(r)^2 d \Omega ^2$.

 We first assume the ambient space $N=I \times G $ where $I=(0,c), 0<c\leq \infty $ and $G$ is equipped with the metric $d \Omega ^2$. $N$\ has warped product metric $ds^2=dr^2+\phi(r)^2 d \Omega ^2$. Let $p$ corresponds the point when $r \rightarrow 0$. (If $G=S^{n-1}$ and $d \Omega ^2$ is standard sphere metric, then the conditions are that $\phi$ is smooth with $\lim_{r \rightarrow 0}\phi = 0$ \ and the Taylor expansion of $\phi$ about $0$ has only odd terms and coefficient for $r^1$ should be 1).

 We consider the case that requires $\phi>0$ to be smooth with $\lim_{r \rightarrow 0}\phi = 0$ and $d \Omega ^2$ \ to be a spherical metric. Then locally use the sphere coordinate $\Theta=\theta_1,...,\theta_{n-1}$ as well as $r$,  local coordinate tangent vector field $\partial_i=\frac{\partial}{\partial \theta_i}, i=1,...,n-1$ and $\partial_n=\partial_r$ are mutually orthogonal. The connection coefficient $\Gamma_{ij}^n$ is proportional to $\phi '/\phi $. Moreover we already have the connection coefficient $\tilde{\Gamma}_{ij}^n$ is $\delta_{ij}/r$ for metric $ds^2=dr^2+r^2 d \Omega ^2$\ because it is a flat metric and whose $r$ level set is obviously a part of a standard sphere, we now can compute $\Gamma_{ij}^n$ $\Gamma_{ij}^n=\tilde{\Gamma}_{ij}^n \cdot \frac{\phi '}{\phi} / \frac{r'}{r} =\frac{\phi '}{\phi}\delta_{ij}$.
 Thus, the second fundamental form of level set $r=\rho_0$ in ambient space $N$ will be $\frac{\phi '}{\phi}I$.
 Using the same argument as in the Lemma1,and now $D$ denotes the covariant derivative in $N$, we have
 \[\sum_{i=1}^{k} \langle D_{\tau_i}W, \tau_i \rangle = kf \frac{\phi '}{\phi} + (f'-f \frac{\phi '}{\phi}) \cdot |(\partial_r)^T|^2\].

If we use \[f(r)=(\phi(r))^{-(k-1)} \int_{\rho_0}^{r}(\phi(\rho))^{k-1} \mathrm d \rho\] and $W=f(r)\partial_r$, then $W$ also satisfies three conditions with little changes in constants. Similar argument can be used to finally derive corresponding lower estimate of area of the minimal surface $M$, $|M|$.

However it should noticed that three remarkable changes will appear in above argument.
\begin{enumerate}
  \item the definition of $N_0$ is not clear when metric is not smooth at $p$, so we just choose constant without interpreting
  \item we have used the fact that $\phi'/\phi \geq 0$ in the hyperbolic case. So we need $\phi'(0) \geq 0 $ and work on the interval $(0,R)$ on which $\phi' \geq 0$
  \item value of the following limit of integral will change \[
\lim_{r \rightarrow 0}\int_{M \cap \partial B_p(r)}\langle W^T, \nu \rangle\].
\end{enumerate}
Regarding the area as the integral of volume form. When the fundamental group of $G$ no longer trivial, we will encounter changes when integrating on equator of $S^{n-1}$ rather than integrating on equator of $G$.

 We may not consider these cases and assume further that $G$ is merely a sphere, i.e. simply connected. Moreover, assume $\phi'(0)>0$ and we must let the radius $r\leq R$ where $\phi' \geq 0$ on $[0,R]$, Then
 \[|W| \rightarrow (r/\phi'(0))^{k-1}\int_{0}^{\rho_0}\phi(\rho)^{k-1} \mathrm d \rho \mathrm  \ \ \mathrm{as} \ \ r \rightarrow 0\]
 Let \[\omega =\int_{0}^{\rho_0}\omega_{k-1}\phi(\rho)^{k-1} \mathrm{d} \rho\]

 \[\lim_{r \rightarrow 0}\int_{M \cap \partial B_p(r)}\langle W^T, \nu \rangle = \lim_{r \rightarrow 0}\int_{M \cap \partial B_p(r)} -\frac{\omega}{\omega_{k-1}}(\frac{r}{\phi'(0)}) ^{-(k-1)}(1+o(r))\]
 Now notice the volume form on $M \cap \partial B_p(r)$ is $\phi ^{k-1}$ times the volume form on $S^{n-1}$, and that $\phi / (\phi'(0)r) \rightarrow 1 $ \ as \  $r\rightarrow 0$. \ The limit of integral is again $\omega$. When the metric is smooth, $\omega$ is just $|N_0|$ as defined by the image of plane in unit ball under exponential map . So we have derived lower bound for minimal surface with boundary on boundary of geodesic ball in the case when $G$ is simply connected sphere.

\end{document}